\documentclass[11pt]{article}
\usepackage{fullpage}
\usepackage{enumerate}
\usepackage{multicol}
\usepackage[usenames,dvipsnames]{color}
\usepackage{hyperref}
\hypersetup{
   hidelinks=true,
   colorlinks=true,
   citecolor=MidnightBlue,
   linkcolor=black
}
\usepackage{comment}
\usepackage{graphicx}
\usepackage{amsmath}
\usepackage{amssymb}
\usepackage{amsthm}
\usepackage{subfigure}
\usepackage{enumitem}
\usepackage{natbib}
\usepackage{dashbox}
\usepackage{fancyvrb}
\usepackage{xcolor}

\definecolor{lightgray}{rgb}{.95,.95,.95}

\newcommand{\minimize}{\mbox{\rm minimize }}

\newcommand{\st}{\mbox{\rm subject to }}

\newcommand{\R}{\mbox{\rm\bf R}}

\long\def\killtext#1{}
\newlength{\lpskip}
\title{MOS: A Mathematical Optimization Service}
\author{
  James Hubert Merrick\footnote{\texttt{jmerrick@alumni.stanford.edu}}, 
  Tom\'{a}s Tinoco De Rubira\footnote{\texttt{ttinoco@alumni.stanford.edu}}}
\date{\today}

\begin{document}

\maketitle

\begin{abstract}
We introduce MOS, a software application designed to facilitate the deployment, integration, management, and analysis of mathematical optimization models. MOS approaches mathematical optimization at a higher level of abstraction than existing optimization modeling systems, enabling its use with all of them. The sole requirement to harness MOS is a simple annotation of the code specifying the formulation of an optimization model. With this, the model becomes accessible to humans through the automatic generation of a user interface, and to machines through an associated API and client libraries. All this is achieved while avoiding the ad hoc code typically required to obtain such features.
\end{abstract}

\section{Introduction}

Whilst known by different names in different settings, mathematical optimization influences billions of dollars in the modern economy, impacting every industry. It consists of the maximization or minimization of some objective function subject to constraints, and is a natural paradigm for solving problems in many fields. In planning applications in operations research and management science, optimization models act as a decision-support tool for a human decision-maker. In other, more operational, settings, they enable decisions to be automated, with example problems including pricing, scheduling, allocation of scarce resources, and routing. In the natural sciences, optimization models capture the physical laws governing the behavior of natural systems. In machine learning, they provide the tools for obtaining parameterized models that best fit a particular data set. \cite{Boyd2004} and \cite{Luenberger2021} provide comprehensive introductions to the theory and applications of mathematical optimization. 

Key tools in mathematical optimization are algebraic modeling systems. Examples of these are \textsc{cvxpy} \citep{diamond2016cvxpy}, \textsc{JuMP} \citep{dunning2017jump}, \textsc{Pyomo} \citep{hart2017pyomo}, \textsc{optmod} \citep{ttinoco2020optmod} and \textsc{GAMS} \citep{bussieck2004}. These tools greatly facilitate the process of constructing and solving optimization models on computers. They allow users to construct optimization problems by writing intuitive mathematical expressions, and can utilize many numerical solvers without the need of custom code that expresses the problem in solver-specific data structures and formats.

From the authors' experience in developing optimization models to support and automate decisions, once a model is formulated using one of the above modeling systems, there is often an additional non-trivial programming exercise required to facilitate a human or application to interact with the model. In the absence of this code, using the model requires familiarity with the model's internal and low-level details, which is seldom documented and user friendly, creating barriers for human users and adding complexity to application pipelines. A custom solution, on the other hand, typically requires time and resources, including dedicated software engineers.

By approaching the modeling problem at a higher level of abstraction than existing tools, and by capturing and standardizing common model properties and structure, MOS provides essential deployment, integration, management and analysis features automatically, removing the need for custom solutions. The sole requirement is a simple annotation of the file containing the model code. This allows a focus, at the development stage, on the core modeling task itself, and at the production stage, on the model usage itself, reducing the barriers to obtaining value from a model.

\cite{guerickeframework} propose a framework for deploying optimization models based on microservice architectures, and highlight a gap between solution methods in literature and solution methods in production environments. MOS also attempts to contribute to the closing of this gap through a proposed concrete and universal model representation, a modular and flexible architecture, and an open-source implementation.\footnote{MOS is available at \url{https://github.com/Fuinn}.}

\section{Model representation}

\begin{figure}
  \centering
  \includegraphics[width=0.4\textwidth]{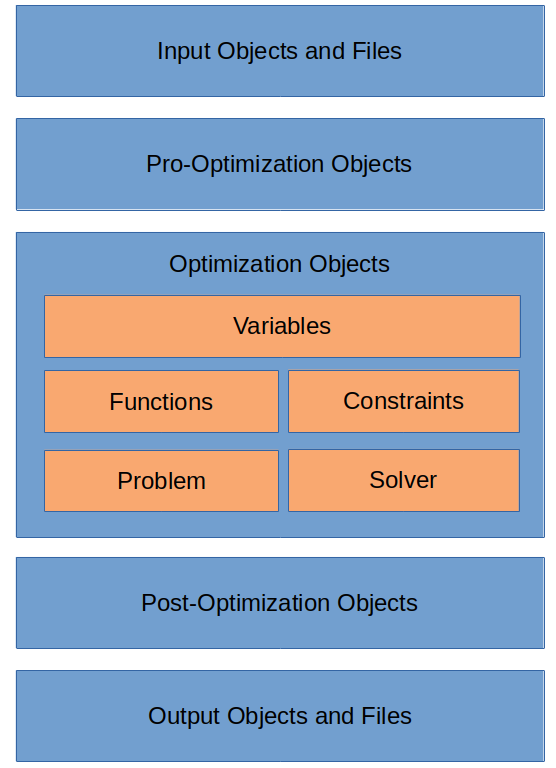}
  \caption{MOS model representation}
  \label{fig:representation}
\end{figure}

\cite{Boyd2004} introduce an optimization problem as being represented by the following:
\begin{equation}\label{structure}
  \begin{array}{rl}
    \minimize & f_0(x) \\
    \st & f_i(x) \leq b_i, \, i=1,\dots,m,
\end{array}
\end{equation}
where $x=(x_1,\dots,x_n)$ is the vector of variables to be optimized, $f_0:\R^n\rightarrow\R$ is the objective function, and functions $f_i:\R^n\rightarrow\R$ together with constants $b_i$ for $i=1,\ldots,m$ define the constraints. MOS considers optimization models as objects that not only consist of optimization problems having variables, functions, and constraints, as in (\ref{structure}), but also of inputs, outputs, and intermediate objects. The intermediate or ``helper objects'' correspond to objects that are created either in a pre-optimization phase, for facilitating the construction of problem variables, functions and constraints, or in a post-optimization phase, for facilitating the construction of outputs. This model representation is helpful for establishing a layer of abstraction that enables the definition and implementation of tools for interacting with, and analyzing models. Figure \ref{fig:representation} illustrates the MOS optimization model representation.

\section{Design}
\label{sec:design}

\begin{figure}
  \centering
  \includegraphics[width=0.618\textwidth]{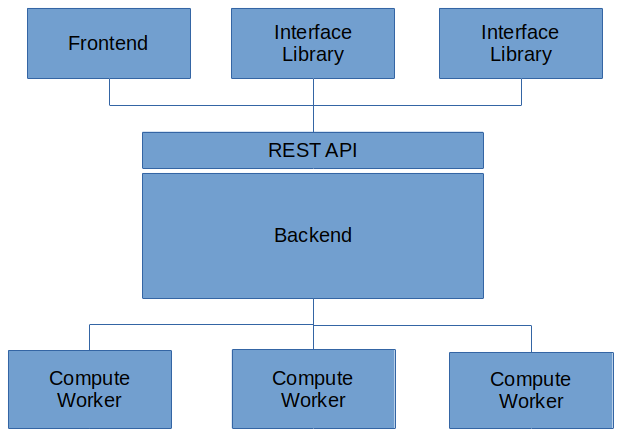}
  \caption{Architecture of MOS}
  \label{fig:architecture}
\end{figure}

Figure \ref{fig:architecture} shows the architecture of MOS, which includes the following components:
\begin{itemize}
\item Backend: Manages model data and access, and provides a REST API for interacting with models.
\item Frontend: Provides a graphical user interface for accessing, using, and analyzing models.
\item Interface libraries: Allow users or other applications to interact with models using popular programming languages and integrate them with application pipelines.
\item Compute workers: Run models locally or distributed over the network using modeling-system-specific computational kernels. 
\end{itemize}

MOS accepts an optimization model file, in any of the supported programming languages and modeling systems, with annotations that identify the components of the model, described in the previous section. Appendix \ref{app:annotate} shows an extract from an annotated model file, illustrating the nature of the annotation, requiring labeling of different sections of the code with certain structural keywords marking out common model features, preceded by a language-specific tag. The degree of structure required by MOS in annotation is variable and optional, with more structure enabling more features in the interface. Appendix \ref{app:api_example} shows an example use of the Python interface library. Appendix \ref{app:s} shows screenshots of the MOS user interface.

\section{Benefits}

MOS enables an appropriately annotated optimization model to have the following features:

 \begin{itemize}

\item \emph{{Deployable:}}
Readily deployable in the cloud or on an organization's own servers.

 \item \emph{API Access:}
Through automatic generation of an API through which an optimization
model may be called, MOS facilitates integration of mathematical optimization
into an organization's data flows and software applications.
For example, MOS {may be used as a microservice by enterprise
  applications to integrate optimization capabilities}.
%
This integrability also enables a modular type of development,
suitable for scaling as
a an optimization model is updated,
without requiring a re-write of supporting custom infrastructure code.

\item \emph{A browsable, graphical, intuitive, representation of a
    model, associated data, and results:} The model structure, documentation, input data, and
  model results, may be browsed through,
  increasing model transparency {and understanding}. Model data and assumptions may also be changed
  through the user interface, and the model solved by clicking a
  button instead of doing so from the command line. 

\item \emph{Leverage of existing optimization technologies and investments:}
  MOS works with a range of domain-specific languages designed for
  optimization such as \textsc{cvxpy}, \textsc{JuMP}, and
\textsc{gams}, not tying an organization to one specific approach as
needs evolve, and additionally leveraging an organization's previous investments in
optimization model code development.

\end{itemize}  


While a well designed custom interface may provide many of these
features, and they are commonly deployed by organizations,
MOS provides these features automatically, while
potentially providing a platform for further benefits. MOS is readily
extendible to incorporate logs of model history, and automated
analysis of model solutions. 



\section{Conclusion}

Mathematical optimization is used
across all industries, and the number of applications are growing as
more digital data is available and the cost of computing declines. {MOS} allows the avoidance of expensive development costs associated
with supporting custom infrastructure code around an optimization
model. It does this by approaching the modelling problem at a higher level of
abstraction than existing tools.



\bibliography{mos_paper}
\bibliographystyle{apalike}

\appendix
\newpage
\section{Extract of annotated model file}
\label{app:annotate}



{\small
\begin{Verbatim}[commandchars=\\\{\}]
\textcolor{ForestGreen}{#@ Constraint: P_limits}
\textcolor{ForestGreen}{#@ Description: Generator active power limits}
P_limits = []
for gen in network.generators:
    P_limits.extend([P[gen.index] >= gen.P_min, P[gen.index] <= gen.P_max])

\textcolor{ForestGreen}{#@ Objective: gen_cost_obj}
gen_cost_obj = optmod.minimize(gen_cost)

\textcolor{ForestGreen}{#@ Problem: problem}
problem = optmod.Problem(gen_cost_obj,
                         constraints=power_balance+angle_ref+P_limits)

\textcolor{ForestGreen}{#@ Solver: solver}
solver = optalg.opt_solver.OptSolverINLP()
solver.set_parameters({`feastol': feastol, `maxiter': 100})

\textcolor{ForestGreen}{#@ Execution: info}
info = problem.solve(solver)

\textcolor{ForestGreen}{#@ Output Object: output_obj}
output_obj = list(info.values())

\end{Verbatim}
}

\section{Example use of MOS Python Interface}
\label{app:api_example}

\begin{verbatim}
from mos.interface import Interface

interface = Interface(url, token)

model = interface.get_model_with_name(`DCOPF Model')

model.set_interface_object(`feastol', 1e-3)
model.set_interface_file(case, `ieee14.m')


model.show_recipe()
model.show_components()
model.run()

model.get_status()
model.get_execution_log()

\end{verbatim}

\section{MOS user interface screenshots}
\label{app:s}
\begin{figure}[!htb]
  \centering
  \dbox{
  \includegraphics[width=0.95\textwidth]{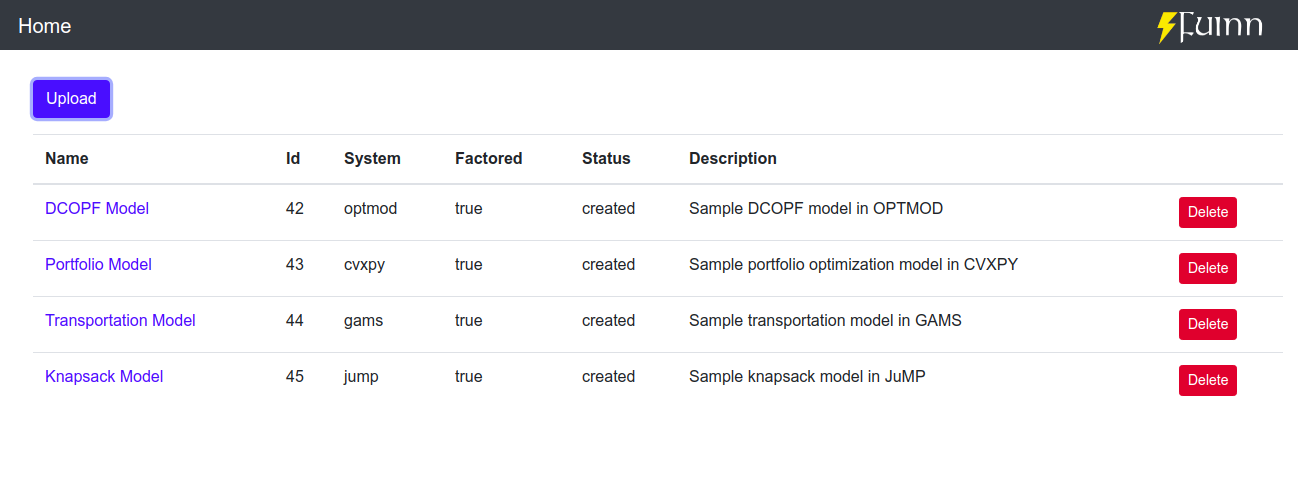}}
  \caption{{MOS} user interface: {list of models available}}
  \label{fig:intro}
\end{figure}

\begin{figure}[!htb]
  \centering
  \dbox{
  \includegraphics[width=0.95\textwidth]{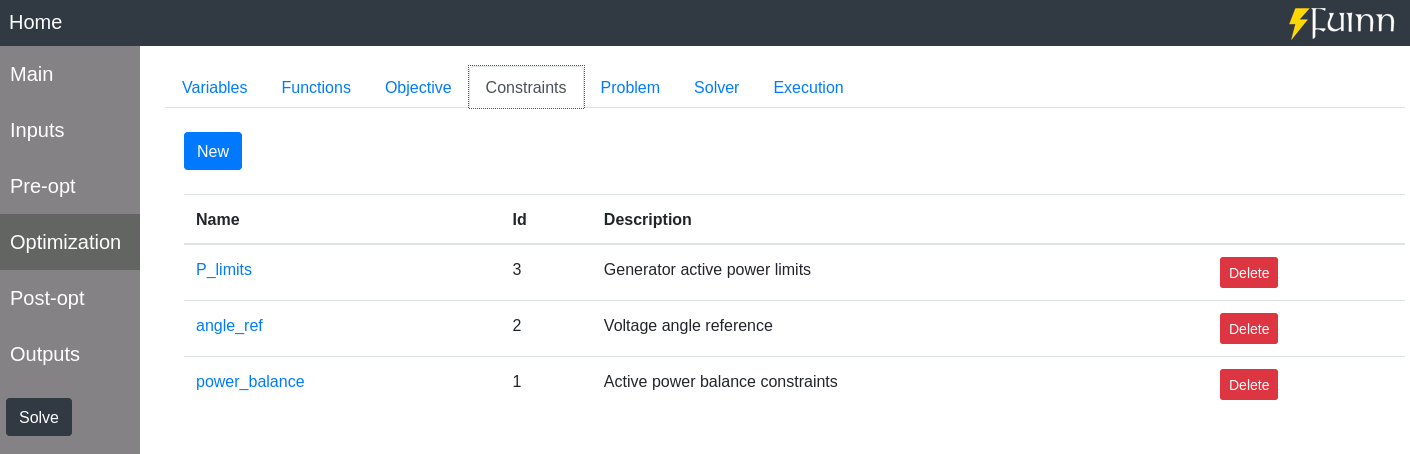}}
  \caption{{MOS} user interface: browsing through model constraints}
  \label{fig:constraints}
\end{figure}

\end{document}